\documentclass[a4paper,10pt]{article}
\usepackage{amsmath,amsthm,amssymb,amsfonts,epsfig,epstopdf,titling,url,array,enumerate,cite,mathrsfs,upgreek,authblk,scrextend, graphicx, float}
\usepackage[bindingoffset=0.2in,left=1in,right=1in,top=1in,bottom=1in,footskip=.25in]{geometry}
\theoremstyle{plain}
\newtheorem{thm}{Theorem}[section]
\providecommand{\keywords}[1]{\begin{addmargin}[28pt]{28pt}\noindent\textbf{Keywords:} #1 \end{addmargin}}
\newtheorem{lma}[thm]{Lemma}

\newtheorem{ppn}[thm]{Proposition}
\theoremstyle{definition}
\newtheorem{dfn}[thm]{Definition}
\newtheorem{eg}[thm]{Example}
\newtheorem{rem}[thm]{\textit{Remark}}
\providecommand{\ams}[1]{\begin{addmargin}[28pt]{28pt}\noindent\textbf{Mathematics Subject Classification:} #1\end{addmargin}}

\title{\textbf{Fuzzy strong $\phi$-b-normed linear space for fuzzy bounded linear operators}}
\author[1]{Abhishikta Das}
\author[2]{T. Bag\footnote{corresponding author}}
\affil[1,2]{Department of Mathematics, Siksha-Bhavana,	\authorcr	Visva-Bharati, Santiniketan-731235, West Bengal, India  
	\authorcr 	E-mail\textsuperscript{1}: abhishikta.math@gmail.com
	\authorcr E-mail\textsuperscript{2}: tarapadavb@gmail.com}
\affil[1]{Orcid Id: https://orcid.org/0000-0002-2860-424X}
\affil[2]{Orcid Id: https://orcid.org/0000-0002-8834-7097}

\date{}

\begin{document}
	\maketitle
	\begin{abstract}
		
		\noindent
	In this paper, concept of fuzzy continuous operator, fuzzy bounded linear operator are introduced in fuzzy strong $\phi$-b-normed linear spaces and their relations are studied. Idea of operator fuzzy norm is developed and  completeness of BF(X,Y) is established.
	\end{abstract}
	\ams{54A40, 03E72, 46A19, 15A03}.
	\section{Introduction}
	The problem of defining fuzzy norm on a linear space was first initiated by Katsaras\cite{8} and afterwards Felbin\cite{9}. Cheng and Mordeson\cite{10}, came up with their definition of fuzzy norm approaching from different perspective. Following the definition of fuzzy norm introduced by Cheng and Mordeson, Bag and Samanta\cite{3} defined a fuzzy norm with a view to establish a complete decomposition of fuzzy norm into crisp norms. However, for doing so, they had to restrict the underlying t-norm in the triangle inequality of fuzzy norm to be the `min'. Next they attempted to deal with the problem by considering general t-norm and established many fundamental results of functional analysis in fuzzy settings\cite{4,5,6,7}. \\
	Following the concept of Bag and Samanta type fuzzy normed linear space, several authors developed the concept of generalized fuzzy normed linear spaces viz. fuzzy cone normed linear space\cite{30}, G-fuzzy normed linear space\cite{31}, etc. In our earlier papers\cite{1,2}, concept of   fuzzy strong $\phi$-b-normed linear space is introduced and developed many basic results in such spaces. In Bag and Samanta type fuzzy normed linear space, scalar multiplication is given by $ N ( cx, t ) = N ( x, \frac{t}{|c|} ) $. But in fuzzy strong $\phi$-b-normed linear space, scalar multiplication is given by $ N ( cx, t ) = N ( x, \frac{t}{\phi (c) } ) $ where $ \phi $ is a real valued function satisfying some properties. \\
	In this paper, we have been able to proceed further. The concept of fuzzy bounded linear operators, fuzzy continuous operators,  operator norm for fuzzy bounded linear operators and spaces of fuzzy bounded linear operators(denoted by BF(X,Y)) are introduced in fuzzy strong $\phi$-b-normed linear spaces. Lastly completeness of BF(X,Y) is  proved. \\
	The organization of the article is in the  following. \\
	Section 2 contains some preliminary results. In Section 3,  fuzzy bounded linear operators, fuzzy continuous operators are  defined and some related results are studied.  Section 4 consists of the study on operator fuzzy norm on   fuzzy strong $\phi$-b-normed linear spaces. 
	%
	%
	\section{Preliminaries}
	We provide some basic definitions and results which are used in this paper.
	\begin{dfn} \cite{3}
		A binary operation $ * : [0 , 1] \times [0 , 1] \rightarrow [0 , 1] $ is called a  $t$-norm if it satisfies the following conditions:
		\begin{enumerate}[(i)]
			\item  $ * $ is associative and commutative; 
			\item $ \alpha ~ * ~ 1 = \alpha ~ ~ \forall \alpha \in [0 , 1] $;  
			\item $ \alpha * \gamma \leq  \beta * \delta ~ $  whenever $ \alpha \leq  \beta  $ and $  \gamma \leq   \delta  ~  ~\forall \alpha, \beta,
			\gamma, \delta \in [0 , 1] $.
		\end{enumerate}
		If $ * $ is continuous then it is called continuous $t$-norm.
	\end{dfn}
	The following are examples of some t-norms. 
	\begin{enumerate}[(i)]
		\item  Standard intersection: $ \alpha *   \beta = \min \{ \alpha ,   \beta \} $.
		\item   Algebraic product: $ \alpha *   \beta = \alpha    \beta $.
		\item Bounded difference: $ \alpha *   \beta = \max \{ 0,  \alpha +  \beta -1 \} $. 
		\item Drastic intersection: $ \alpha *   \beta =  \begin{cases} 
			\alpha ~~~~ if ~ \beta = 1 \\
			\beta ~~~~ if ~ \alpha = 1 \\
			0 ~~~~ otherwise 
		\end{cases} $
	\end{enumerate}
	\begin{dfn} \cite{5}
		Let $ X $ be a linear space over a field $ \mathbb{F} $. A  fuzzy subset $ N $ of $ X \times \mathbb{R}$ 
		is called fuzzy norm on $  X $ if  $ \forall x, y \in X $ the following conditions hold:
		\begin{enumerate}[(NI)] 
			\item $ \forall t \in {\mathbb{R}}  $ with $ t \leq 0, ~ N(x,t)=0 $; 
			\item $ (\forall t \in {\mathbb{R}}, ~ t>0, ~ N(x,t) = 1)   \iff  x = \theta $;
			\item $ \forall t \in {\mathbb{R}}, $ and  $ c \in {\mathbb{R}}, ~ t>0, ~ N(cx,t) = N(x, \frac{t}{|c|}) $;
			\item $ \forall s, t \in {\mathbb{R}}, ~ N(x+y, s+ t) \geq N(x,s) * N(y,t) $;
			\item $ N(x, \cdot) $ is a non-decreasing function of $ t $ and $ \underset{t \rightarrow \infty }{lim} N(x,t) = 1 $. 
		\end{enumerate}
		Then the pair  $ ( X, N ) $ is called  fuzzy normed linear space. 	
	\end{dfn}
	\begin{rem} 
		Bag and Samanta \cite{3} assumed that,
		$$ ~~ \textbf{(NVI)}: ~ N(x, t ) > 0 ~ \implies ~ x = \theta. ~ $$ 
	\end{rem}
	\begin{dfn} \cite{3}
		Let $ ( X, N ) $  be a fuzzy normed linear space.
		\begin{enumerate}[(i)]
			\item A sequence $ \{ x_n \} $ is said to be convergent if $ \exists x \in X $ 
			such that $ \underset{n \rightarrow \infty }{lim} N(x_n - x, t) = 1, ~ \forall t > 0 $. Then  $ x $ is called the limit of the sequence $ \{ x_n \} $
			and denoted by $ lim x_n $.
			\item A sequence $ \{ x_n \} $ in a fuzzy normed linear space $ ( X, N ) $ is said to be Cauchy if \\
			$ \underset{n \rightarrow \infty }{lim} N(x_{n+p} - x_n, t) = 1, ~ \forall t > 0 $ and $ p = 1,2, \cdots $.
			\item  $ A \subseteq  X $ is said to be a closed set if for any sequence $ \{ x_n \} $ in $ A $ converging to $ x \in X $ implies $ x \in A $. 
			\item $ A \subseteq  X $ is said to be the closure of $ A $, denoted by $ \bar{A} $ if for any $x \in \bar{A} $, there is a sequence
			$ \{ x_n \} \subseteq A $ such that $ \{ x_n \} $ converges to  $ x $. 
			\item $ A \subseteq  X $ is said to be compact if any sequence $ \{ x_n \} \subseteq A $  has a subsequence converging to an element in $ A $.
		\end{enumerate}
	\end{dfn}
	\begin{rem}
		For a t-norm, $ *: [0,1] \times [0,1] \rightarrow [0,1] $, Bag and Samanta\cite{7} assumed
		$$ \textbf{(T1)} ~~ \text{for}~ ~ \alpha > 0, ~ \alpha * \alpha > 0. $$
	\end{rem}
	\begin{dfn}\cite{7}
		Let $ ( X, N ) $ be a fuzzy normed linear space and $ \alpha \in (0, 1 ) $.
		\begin{enumerate}[(i)]
			\item A sequence $ \{ x_n \} $ in $ X $ is said to be $\alpha $-fuzzy convergent if $ \exists x \in X $ such that 
			$$ \underset{n\rightarrow \infty}{lim}  \bigwedge \{ t > 0 : N (x_n-x, t) > 1 -\alpha  \} = 0. $$
			If $ \{ x_n \} $ is  $\alpha $-fuzzy convergent for each $ \alpha \in (0, 1 ) $ then $ \{ x_n \} $ is called level fuzzy (l-fuzzy) convergent.
			\item A sequence $ \{ x_n \} $ in $ X $ is said to be $\alpha $-fuzzy Cauchy if $ \exists x \in X $ such that 
			$$ \underset{n,m\rightarrow \infty}{lim}  \bigwedge \{ t > 0 : N (x_n-x_m, t) > 1 -\alpha  \} = 0. $$
			If $ \{ x_n \} $ is  $\alpha $-fuzzy Cauchy for each $ \alpha \in (0, 1 ) $ then $ \{ x_n \} $ is called level fuzzy (l-fuzzy) Cauchy.
			\item A subset $ F  $ of $ X $ is said to be $\alpha $-fuzzy complete if every $\alpha $-fuzzy Cauchy sequence is $\alpha $-fuzzy convergent
			to some point in $ F $. \\
			If $F $ is  $\alpha $-fuzzy complete for each $ \alpha \in (0, 1 ) $ then $F $ is called level fuzzy (l-fuzzy) complete.
			\item A subset $ F  $ of $ X $ is said to be $\alpha $-fuzzy closed if for any sequence $ \{ x_n \} $ in $ F $ is  $\alpha $-fuzzy convergent
			to $ x \in X $ implies $ x \in  F $. \\
			If $F $ is  $\alpha $-fuzzy closed for each $ \alpha \in (0, 1 ) $ then $F $ is called level fuzzy (l-fuzzy) closed.
			\item A set $ S $ is said to be $\alpha $-fuzzy closure of  $ F  \subset X $ if for each $ x \in S $, there exist a $\alpha $-fuzzy convergent 
			sequence converging to $ x $. \\
			It is denoted by   $\alpha $-$\bar{ F } $. 
			\item A subset $ F  $ of $ X $ is said to be $\alpha $-fuzzy compact if for any sequence $ \{ x_n \} $ in $ F $, there exists  a subsequence of 
			$ \{ x_n \} $ which is $\alpha $-fuzzy convergent to some point in $ F $. \\
			If $F $ is  $\alpha $-fuzzy compact for each $ \alpha \in (0, 1 ) $ then $F $ is called level fuzzy (l-fuzzy) compact.
			\item A subset $ F  $ of $ X $ is said to be $\alpha $-fuzzy  bounded if there exists $  t(\alpha) > 0 $ such that\\ $ N(x,t) > 1 -\alpha, ~
			\forall x \in F $. \\
			If $F $ is  $\alpha $-fuzzy  bounded for each $ \alpha \in (0, 1 ) $ then $F $ is called level fuzzy (l-fuzzy)  bounded.
		\end{enumerate}
	\end{dfn}
\begin{rem}
	In \cite{6}, $ \bigwedge \{ t > 0 : N_1 (x, t ) \geq \alpha \}, ~ \alpha \in (0,1) $ is denoted by $ d_\alpha $. By (NVI), it follows that for $ x \neq \theta, ~ d_\alpha > 0 ~~  \forall \alpha \in (0,1)  $. 
\end{rem}
	\begin{dfn}\cite{1}
		Let  $ \phi $ be a function defined on  $ {\mathbb{R}} $ to $ {\mathbb{R}} $ with the following properties 
		\begin{enumerate}[($\phi$1)]
			\item $ \phi (-t) = \phi (t), ~ \forall t \in {\mathbb{R}}  $; 
			\item $ \phi (1) = 1 $;
			\item $ \phi $ is strictly increasing and continuous on $ (0, \infty) $;
			\item $ \underset{\alpha \rightarrow 0 }{lim} \phi(\alpha) = 0 $ and $ \underset{\alpha \rightarrow \infty }{lim} \phi(\alpha) = \infty $.
		\end{enumerate}
	\end{dfn}
	The following are examples of such functions.
	\begin{enumerate}[(i)]
		\item $ \phi(\alpha) = | \alpha |,  ~ \forall \alpha \in {\mathbb{R}}  $
		\item $ \phi(\alpha) = | \alpha |^p, ~\forall \alpha \in {\mathbb{R}}, ~p \in \mathbb{R}^+  $
		\item $ \phi(\alpha) = \frac{2 \alpha^{2n}}{|\alpha| + 1}, ~ ~\forall \alpha \in {\mathbb{R}}, ~ n \in \mathbb{N} $
	\end{enumerate}
	\begin{dfn}\label{dfn1}\cite{1}
		Let $X$ be a linear space over the field $ \mathbb{R} $ and $ K \geq 1 $  be a given real number. A  fuzzy subset $ N $ of $ X \times \mathbb{R}$ 
		is called fuzzy strong $\phi$-b-norm on $  X $ if  for all $ x, y \in X $ the following conditions hold:
		\begin{enumerate}[(bN1)] 
			\item $ \forall t \in {\mathbb{R}}  $ with $ t \leq 0, ~ N(x,t)=0 $; 
			\item $ (\forall t \in {\mathbb{R}}, ~ t>0, ~ N(x,t) = 1)  $ iff $ x = \theta $;
			\item $ \forall t \in {\mathbb{R}}, ~ t>0, ~ N(cx,t) = N(x, \frac{t}{\phi(c)} ) $ if $\phi(c) \neq 0 $;
			\item $ \forall s, t \in {\mathbb{R}}, ~ N(x+y, s+ Kt) \geq N(x,s) * N(y,t) $;
			\item $ N(x, \cdot) $ is a non-decreasing function of $ t $ and $ \underset{t \rightarrow \infty }{lim} N(x,t) = 1 $. 
		\end{enumerate}
		Then $ ( X, N, \phi, K, * ) $ is called  fuzzy strong $\phi$-b-normed linear space. 	
	\end{dfn}

	\begin{lma} \cite{2} \label{pre lma1}
		Let $  ( X, N, \phi, K, * ) $ be a fuzzy strong $\phi$-b normed linear space satisfying (N6) and the underlying $t$-norm $*$ is continuous at 
		$ (1,1) $. If $ \{ x_1, x_2, \cdots x_n \} $ is  a linearly independent set of vectors in $ X $, then for each 
		$ \alpha \in (0,1), ~ \exists c_\alpha > 0 $ such that  for any set of scalars  $ \{ \beta_1, \beta_2, \cdots \beta_n \} $ with
		$ \sum _{i=1} ^n |\beta_i| \neq 0 $, 
		\begin{equation*}
			\bigwedge \{ t > 0 : N( \beta_1x_1 +  \beta_2x_2 + \cdots + \beta_nx_n, Kt  )\geq 1- \alpha  \} \geq  
			\frac{c_\alpha} {\phi( \frac{1}{\sum _{i=1} ^n |\beta_i|})}.
		\end{equation*}  
	\end{lma}
		%
		%

\section{Fuzzy bounded linear operator}
In this section   definition of  fuzzy bounded linear operator and fuzzy continuous operator in fuzzy strong $\phi$-b-normed linear space are given and their relation is studied.
\begin{dfn}
	Let $ T : (X, N_1, \phi, K, *_1 ) \rightarrow  (Y, N_2, \phi, K, *_2) $ be a linear operator where $ X $ and $Y$ are fuzzy strong $\phi$-b-normed linear spaces. $T$ is said to be
	fuzzy bounded if for each $ \alpha \in (0,1), ~~ \exists M_\alpha > 0 $ such that 
	\begin{equation}  \label{equn 1}
		N_1 (x, \frac{t}{M_\alpha})  \geq 1 - \alpha \implies N_2(Tx, Ks) \geq \alpha ~~ \forall s > t ~ \text{and} ~ \forall t> 0. 
	\end{equation}
\end{dfn}
\begin{ppn} \label{ppn 1}
	Let $ T : X \rightarrow  Y $  be a linear operator where $ X $ and $Y$ are fuzzy strong $\phi$-b-normed linear spaces. If $T$ is  fuzzy bounded, then the relation
	(\ref{equn 1}) is equivalent to the relation 
	\begin{equation}  \label{equn 2}
	\bigwedge \{ t>0 : N_2 (Tx, K t ) \geq \alpha \} \leq M_\alpha 	\bigwedge  \{ t>0 : N_1 (x, t ) \geq 1 - \alpha \} ~ ~ \forall x \in X. 	
	\end{equation}
\end{ppn}	
\begin{proof}
First we show that $ (\ref{equn 1}) \implies (\ref{equn 2}) $. 
\begin{align*}
\text{Let}  ~ & r > 	M_\alpha 	\bigwedge \{ t>0 : N_1 (x, t ) \geq 1 - \alpha \}  \\
\implies & \dfrac{r}{M_\alpha} > 	\bigwedge \{ t>0 : N_1 (x, t ) \geq 1 - \alpha \}  \\ 
\implies & \exists \dfrac{r'}{M_\alpha} ~ \text{such that} ~ \dfrac{r}{M_\alpha} > \dfrac{r'}{M_\alpha}  > \bigwedge \{ t>0 : N_1 (x, t ) \geq 1 - \alpha \}  \\	
\implies & \exists \dfrac{r'}{M_\alpha} ~ \text{such that} ~ N_1 ( x, \dfrac{r'}{M_\alpha} ) \geq 1 - \alpha \\
\implies &  N_2 ( Tx, K r ) \geq \alpha ~~~~~              ( ~ \text{by (\ref{equn 1})} ~ ) \\
\implies & 	\bigwedge \{ t>0 : N_2 (Tx, K t ) \geq \alpha \} \leq r \\
\implies & 	\bigwedge \{ t>0 : N_2 (Tx, K t ) \geq \alpha \} \leq  M_\alpha 	\bigwedge  \{ t>0 : N_1 (x, t ) \geq 1 - \alpha \} ~ ~ \forall x \in X.  
\end{align*}
So, $ (\ref{equn 1}) \implies (\ref{equn 2}) $. \\
Now we prove $ (\ref{equn 2}) \implies (\ref{equn 1}) $. \\
Assume that $ 	N_1 (x, \frac{t}{M_\alpha} )  \geq 1 - \alpha $. So, 
\begin{align*}
&	\bigwedge  \{ r > 0 : N_1 (x, r ) \geq 1 - \alpha \} \leq   \frac{t}{M_\alpha} \\
\implies &  M_\alpha  \bigwedge  \{ r > 0 : N_1 (x, r ) \geq 1 - \alpha \} \leq   t \\
\implies & \bigwedge  \{ r > 0 : N_2 ( Tx, K r ) \geq \alpha \} \leq   t     ~~~~~  ( ~ \text{by (\ref{equn 2}) } )   \\ 
\implies & \text{for any }  ~ s > t,  ~~ \bigwedge  \{ r > 0 : N_2 ( Tx, K r ) \geq \alpha \}  < s \\
\implies & N_2 ( Tx, K s ) \geq \alpha.
\end{align*}
Hence, $ (\ref{equn 2}) \implies (\ref{equn 1}) $.
\end{proof}	
\begin{rem}
	We denote the collection of all linear operators defined from a fuzzy strong $\phi$-b-normed linear space  $ (X, N_1, \phi, K, *_1 ) $ to another  $  (Y, N_2, \phi, K, *_2 ) $
	by $ L(X, Y) $ and for fuzzy bounded linear operators we denote  the collection by  $ BF (X, Y) $.
\end{rem}	
\begin{lma} \label{lma 1}
	Let $ (X, N, \phi, K, * ) $ be a fuzzy strong $\phi$-b-normed linear space and the underlying t-norm $ * $ be continuous at $ ( 1, 1 ) $. Then for each $ \alpha \in (0, 1) $,
	$ \exists \beta \geq \alpha $ such that 
	$$ \bigwedge \{ s+t > 0 : N ( x+y, K (s+t) ) \geq \alpha \} \leq   \bigwedge \{ s > 0 : N ( x, s ) \geq \beta \} +   \bigwedge \{ t > 0 : N ( y, t ) \geq \beta \}, 
	\forall x, y \in X $$ 
\end{lma}	
\begin{proof}
	 Since $ * $ is continuous at $ ( 1, 1 ) $, for each  $ \alpha \in (0, 1) $, we can find  $ \beta \in (0, 1) $ such that $ \beta * \beta \geq \alpha $. \\
	 Again $ \beta \geq \beta * \beta \geq \alpha \implies \beta \geq \alpha $. Now,
	 \begin{align*}
	 	 & \bigwedge \{ s > 0 : N ( x, s ) \geq \beta \} +   \bigwedge \{ t > 0 : N ( y, t ) \geq \beta \} \\
	 	 = & \bigwedge \{ s+ t > 0 :  N ( x, s ) \geq \beta, ~ N ( y, t ) \geq \beta \} \\
	 	 \geq &  \bigwedge \{ s+ t > 0 :  N ( x, s ) * N ( y, t ) \geq \beta * \beta \} \\
	 	 \geq &  \bigwedge \{ s+ t > 0 : N ( x+y, s + K t ) \geq \alpha \}
	 \end{align*}
 Hence $ \forall x, y \in X $,
 \begin{equation} \label{equn 3}
 	\bigwedge \{ s+ t > 0 : N ( x+y, s + K t ) \geq \alpha \} \leq \bigwedge \{ s > 0 : N ( x, s ) \geq \beta \} +   \bigwedge \{ t > 0 : N ( y, t ) \geq \beta \} 
 \end{equation}
Now,  $ K \geq 1 \implies  K s \geq s \implies  K s + K t \geq s + K t \implies K (s +t ) \geq s + K t $. So, 
\begin{align*}
& \{ s+ t > 0 : N ( x+y, s + K t ) \geq \alpha \} \subset \{ s+ t > 0 : N ( x+y, K (s +t ) ) \geq \alpha \}	\\
\implies & 	\bigwedge \{ s+ t > 0 : N ( x+y, s + K t ) \geq \alpha \} \geq \bigwedge \{ s+ t > 0 : N ( x+y, K (s +t ) ) \geq \alpha \} 
\end{align*}
From ( \ref{equn 3}) we get, 
$$ \bigwedge \{ s+ t > 0 : N ( x+y, K (s +t ) ) \geq \alpha \} \leq \bigwedge \{ s > 0 : N ( x, s ) \geq \beta \} +   \bigwedge \{ t > 0 : N ( y, t ) \geq \beta \},
 \forall x, y \in X $$
\end{proof}	
\begin{thm} \label{thm 1}
	$ BF (X, Y ) $ (set of all fuzzy bounded linear operators) is  a subspace of $ L (X, Y ) $(set of all  linear operators) where $  (X, N_1, \phi, K,  *_1) $ and 
	$  (Y, N_2, \phi, K, *_2) $ are fuzzy strong $\phi$-b-normed linear spaces and  $ *_2 $ is continuous at $ ( 1, 1 ) $.  
\end{thm}	
\begin{proof}
	We take $ T_1, T_2 \in BF (X, Y ) $. \\
	Now by Lemma (\ref{lma 1}),  for non-zero scalars $ k_1, k_2 $, we have
	 \begin{align*}
	 	\bigwedge \{ s+ t > 0 : N_2 ( (k_1 T_1 + k_2 T_2) x, K (s +t ) ) \geq \alpha \} \leq & \bigwedge \{ s > 0 : N_2 ( k_1 T_1 x, s ) \geq \beta \} + \\
	 	 & \bigwedge \{ t > 0 : N ( k_2 T_2 x, t ) \geq \beta \}, ~ 
	 	 \forall x \in X
	 \end{align*}   
	  where $ \beta $ depends on $ \alpha $ and $ \beta \geq \alpha $. \\
	  Therefore,
	  $$ 	\bigwedge \{ s+ t > 0 : N_2 ( (k_1 T_1 + k_2 T_2) x, K (s +t ) ) \geq \alpha \} \leq \phi ({k_1}) \bigwedge \{ s > 0 : N_2 ( T_1 x, s ) \geq \beta \}  $$
	  \begin{equation} \label{equn 4}
	  ~~~~~~~~~~~~~~~~~~~~~~~~~~~~~~~~~~~~~~~~~~~~~~~~~~~~~~~~~~	+ \phi ({k_2}) \bigwedge \{ t > 0 : N (  T_2 x, t ) \geq \beta \}, ~ 	\forall x \in X
	  \end{equation}
	  Since  $ T_1, T_2 $ are fuzzy bounded, $ \exists ~  M_{\beta(\alpha)} ^1,  M_{\beta(\alpha)} ^2 > 0 $ such that 
	  $$ \bigwedge \{ s > 0 : N_2 (  T_1 x, s ) \geq \beta \}  \leq  M_{\beta(\alpha)} ^1 \bigwedge \{ s > 0 : N_1 (   x, s ) \geq 1 - \beta \}, ~ 	\forall x \in X  $$ 
	  and 
	    $$ \bigwedge \{ t > 0 : N_2 (  T_2 x, t ) \geq \beta \}  \leq  M_{\beta(\alpha)} ^2 \bigwedge \{ t > 0 : N_1 (   x, t ) \geq 1 - \beta \}, ~ 	\forall x \in X.  $$ 
	 Now from (\ref{equn 4}) we have,
	 \begin{align*}
	 		\bigwedge \{ s+ t > 0 : N_2 ( (k_1 T_1 + k_2 T_2) x, K (s +t ) )
	 		& \geq \alpha \}   \leq  \phi ({k_1})  M_{\beta(\alpha)} ^1
	 		 \bigwedge \{ s > 0 : N_1 (  x, s ) \geq 1 - \beta \} \\
	 		 & + \phi ({k_2}) M_{\beta(\alpha)} ^2 \bigwedge \{ t > 0 : N_1 (  x, t ) \geq 1 - \beta \}, ~ 	\forall x \in X
	 \end{align*} 
 Let $ M_\alpha = \phi ({k_1})  M_{\beta(\alpha)} ^1 +   \phi ({k_2}) M_{\beta(\alpha)} ^2 $. Then from above we get, 
 \begin{equation} \label{equn 5}
 	\bigwedge \{ s+ t > 0 : N_2 ( (k_1 T_1 + k_2 T_2) x, K (s +t ) ) \geq \alpha \} \leq M_\alpha \bigwedge \{ t > 0 : N_1 (  x, t ) \geq 1 - \beta \}, ~ 	\forall x \in X
 \end{equation}
Since $\beta \geq \alpha $, so $ 1 - \alpha \geq 1 -\beta $. Thus, 
\begin{align*}
	& \{ t > 0 : N_1 (  x, t ) \geq 1 - \alpha \} \subset \{ t > 0 : N_1 (  x, t ) \geq 1 - \beta \} \\
\implies & \bigwedge \{ t > 0 : N_1 (  x, t ) \geq 1 - \alpha \} \geq \bigwedge \{ t > 0 : N_1 (  x, t ) \geq 1 - \beta \} 	
\end{align*}
So from (\ref{equn 5}) we have, 
\begin{align*}
	&	\bigwedge \{ s+ t > 0 : N_2 ( (k_1 T_1 + k_2 T_2) x, K (s +t ) ) \geq \alpha \} \leq  M_\alpha 
		 \bigwedge \{ t > 0 : N_1 (  x, t ) \geq 1 - \alpha \}, ~ 	\forall x \in X \\
		 \implies & k_1 T_1 + k_2 T_2 \in BF (X, Y ) 
\end{align*}
Therefore $ BF (X, Y ) $ is  a subspace of $ L (X, Y ) $.
\end{proof}	
\begin{dfn}
	 An operator $ T : (X, N_1, \phi, K, *_1) \rightarrow  (Y, N_2, \phi, K,   *_2) $ is said to be fuzzy continuous  at $ x \in X $ if for every sequence $ \{ x_n \} $ 
	 in $ X $ with $ x_n \rightarrow x $ implies $ T x_n \rightarrow T x $.\\
	 That is $ \underset{n \rightarrow \infty}{lim} N(x_n - x, t ) = 1, ~ \forall t> 0 \implies  \underset{n \rightarrow \infty}{lim} N( T x_n - T x, t ) = 1, ~ \forall t> 0 $.
\end{dfn}	
\begin{thm} \label{thm 2}
	Let $ T : (X, N_1, \phi, K, *_1) \rightarrow  (Y, N_2, \phi, K,   *_2) $ be a linear operator where $ X $ and $Y$ are fuzzy strong $\phi$-b-normed linear spaces. If $T$ is 
	fuzzy continuous  at a point  $ x_0 \in X $, then $T$ is  	fuzzy continuous  everywhere in $ X $.
\end{thm}	
\begin{proof}
	Proof is straightforward(same as in \cite{6}).
\end{proof}	
\begin{thm}  \label{thm 3}
		Let $ T : (X, N_1, \phi, K,   *_1) \rightarrow  (Y, N_2, \phi, K,   *_2) $ be a linear operator where $ X $ and $ Y $ are fuzzy strong $\phi$-b-normed linear spaces. 
		If $T$ is fuzzy bounded then it is fuzzy continuous but not conversely. 
\end{thm}	
\begin{proof}
	First we suppose that $ T $ is fuzzy bounded. So for each $ \alpha \in (0,1) $, $ \exists M_\alpha $ 
	such that 
	$$ \bigwedge \{ t > 0 : N_2 (  Tx, K t ) \geq \alpha \} \leq M_\alpha \bigwedge \{ t > 0 : N_1 (  x, t ) \geq 1 - \alpha \}, ~ 	\forall x \in X $$
	Let $ \{ x_n\} $ be a sequence in $X $ such that $  x_n \rightarrow x $. Thus, 
	$$  \underset{n \rightarrow \infty}{lim} N_1 (x_n - x, t ) = 1, ~ \forall t> 0.  $$
	Let $ \epsilon > 0 $ be given. So for each $ \alpha \in (0,1) $, $ \exists $ a positive integer $ N (\alpha, \epsilon ) $ such that 
\begin{align*}
	&  N_1 (x_n - x, \dfrac{\epsilon}{2K M_\alpha} ) > 1 - \alpha, ~~~ \forall n \geq N (\alpha, \epsilon ) \\
	 \implies &  \bigwedge \{ t > 0 : N_1 (  x_n - x, t ) \geq 1 - \alpha \} \leq \dfrac{\epsilon}{2K M_\alpha}, ~~ \forall n \geq N (\alpha, \epsilon ),
	 ~ \forall \alpha \in (0,1)  \\
	 \implies & M_\alpha \bigwedge \{ t > 0 : N_1 (  x_n - x, t ) \geq 1 - \alpha \} \leq \dfrac{\epsilon}{2K } < \dfrac{\epsilon}{K }, ~~ \forall n \geq N (\alpha, \epsilon ), ~ \forall \alpha \in (0,1)  \\ 
	 \implies &  N_2 ( Tx_n - Tx, \epsilon ) \geq \alpha , ~~ \forall n \geq N (\alpha, \epsilon ), ~\forall \alpha \in (0,1)  \\ 
	 \implies & \underset{n \rightarrow \infty}{lim} N_2 ( Tx_n - Tx, \epsilon ) = 1
\end{align*} 
	Since  $ \epsilon > 0 $ is arbitrary, we have $ \underset{n \rightarrow \infty}{lim} N_2 ( Tx_n - Tx, t ) = 1, ~ \forall t > 0 $. Hence $ T $ is fuzzy continuous on $X $. \\
	Converse result may not be true which is justified by the following example.
\end{proof}	
\begin{eg}
	Let us consider a normed linear space $ ( X, ||.||)  $. We define two functions as in the following: 
	$$ N_1 (x, t ) = \begin{cases}
		1 ~~~~ if ~~~ t > ||x|| \\
		\dfrac{1}{2} ~~~~ if ~~ 0 < t < ||x|| \\
		0 ~~~~~ if ~~ t \leq 0
	\end{cases} $$
$$ N_2 (x, t ) = \begin{cases}
	1 - \frac{||x||}{t}~~~~ if ~~~ t \geq ||x|| \\
	0 ~~~~~~~~~~~~ if ~~ t < ||x||
\end{cases} $$
Then $ N_1, N_2 $ are fuzzy norm on $X$ where $ * $ is taken as `min'(please see Observation 1.2 \cite{99}).  \\
So $ N_1 $ and $ N_2 $ are strong fuzzy $\phi$-b-norm on $X$ where $ * $ is `min', $ \phi(\alpha) = |\alpha| ~ ~ \forall \alpha \in 
\mathbb{R} $ and $ K= 1$. \\
Now define an linear operator $ T : ( X, N_1 ) \rightarrow ( X, N_2 ) $ by $ T x = 2x ~ ~ \forall x \in X $. Then it can be shown that $T$ is fuzzy continuous but not fuzzy bounded(please see Example 3.1 \cite{6}).
\end{eg}	
\begin{thm} \label{thm 6}
	Let $	T : X \rightarrow Y $ be a linear operator where $ (X, N_1, \phi, K, *_1) $ and $  (Y, N_2, \phi, K, *_2) $ be two fuzzy strong $\phi$-b-normed linear spaces. If $X$ is finite dimensional then $T$ is fuzzy bounded.
	\begin{proof}
		Let $ dim X = n $ and $ \{ x_1, x_2, \cdots , x_n \}  $ be a basis of $X$. \\
		Choose $ x ( \neq \theta ) \in X $. Then $ x = \beta_1 x_1 + \beta_2 x_2 + \cdots + \beta_n x_n $ for some suitable scalars 
		$ \{ \beta_1, \beta_2, \cdots, \beta_n \} $. So, 
		\begin{align*}
			T x = & ~ T ( \beta_1 x_1 + \beta_2 x_2 + \cdots + \beta_n x_n ) \\
			= & ~ T (\beta_1 x_1) +  T (\beta_2 x_2) +  \cdots +  T (\beta_n x_n)
		\end{align*}
		Let $ s = \frac{1}{\sum_{i=1}^{n} | \beta_i | } $. Clearly, $ s \neq 0 $. \\
		Let us choose $ \alpha \in (0,1) $ arbitrary. Now, 
		\begin{align*}
			& \bigwedge \{ \frac{t}{n} > 0 : N_2 ( s T(\beta_1 x_1), \frac{t}{n}  ) \geq \alpha   \} + 
			\bigwedge \{ \frac{t}{n} > 0 : N_2 ( s T(\beta_2 x_2), \frac{t}{K n}  ) \geq \alpha   \} +  \cdots  \\
			& \hskip 160 pt +	\bigwedge \{ \frac{t}{n} > 0 : N_2 ( s T(\beta_n x_n), \frac{t}{K n}  ) \geq \alpha   \} \\
			= &  \bigwedge \{ \frac{t}{n} + \frac{t}{n} > 0 : N_2 ( s T(\beta_1 x_1), \frac{t}{n}  ) \geq \alpha, ~
			N_2 ( s T(\beta_2 x_2), \frac{t}{K n}  ) \geq \alpha   \} +	\cdots  \\
			& \hskip 160 pt  +	\bigwedge \{ \frac{t}{n} > 0 : N_2 ( s T(\beta_n x_n), \frac{t}{K n}  ) \geq \alpha   \} \\
			\geq & \bigwedge \{ \frac{t}{n} + \frac{t}{n} > 0 : N_2 ( s T(\beta_1 x_1) + s T(\beta_2 x_2) , \frac{t}{n} + K \cdot 
			\frac{t}{K n}   ) \geq \alpha * \alpha \} +	\cdots  \\
			& \hskip 160 pt  +	\bigwedge \{ \frac{t}{n} > 0 : N_2 ( s T(\beta_n x_n), \frac{t}{K n}  ) \geq \alpha   \} \\
			& \cdots \hskip 20 pt  \cdots \\
			\geq & \bigwedge \{ \frac{t}{n} +  \cdots +  \frac{t}{n} > 0 : N_2 ( s T(\beta_1 x_1) +  \cdots 
			s T(\beta_n x_n), \frac{t}{n} +   \cdots +  \frac{t}{n} ) \geq \alpha *  \cdots * \alpha \} \\
			= & \bigwedge \{ t > 0 :  N_2 ( T( \beta_1 x_1 + \beta_2 x_2 + \cdots + \beta_n x_n ), \frac{t}{\phi(s)} )  
			\geq \alpha * \alpha * \cdots * \alpha \}
		\end{align*}
		Hence we obtain, 
		\begin{equation*}
			\bigwedge \{ \frac{t}{n} > 0 : N_2 ( s T(\beta_1 x_1), \frac{t}{n}  ) \geq \alpha   \} +  \cdots  
			+	\bigwedge \{ \frac{t}{n} > 0 : N_2 ( s T(\beta_n x_n), \frac{t}{K n}  ) \geq \alpha   \} \geq 
		\end{equation*}
		\begin{equation}  \label{equn 12}
			\hskip 100 pt \bigwedge \{ t > 0 :  N_2 ( T( \beta_1 x_1  + \cdots + \beta_n x_n ), \frac{t}{\phi(s)} )  
			\geq \alpha * \cdots * \alpha \}
		\end{equation}
		Let, 
		$ N_\alpha = \max \{ \bigwedge \{ \frac{t}{n} > 0 : N_2 ( s T(\beta_1 x_1), \frac{t}{n}  ) \geq \alpha   \},  \cdots,	\bigwedge \{ \frac{t}{n} > 0 : N_2 ( s T(\beta_n x_n), \frac{t}{K n}  ) \geq \alpha   \}  \} $.  \\
		Now from (\ref{equn 12}), we get 
		$$ n N_\alpha \geq \bigwedge \{ t > 0 :  N_2 ( Tx, \frac{t}{\phi(s)} )  \geq \alpha * \cdots * \alpha \} $$
		which implies 
		\begin{equation}  \label{equn 13}
			\frac{	n N_\alpha}{\phi(s)} \geq \bigwedge \{ t > 0 :  N_2 ( Tx, t )  \geq \alpha * \cdots * \alpha \}
		\end{equation}
		From Lemma(\ref{pre lma1}), $ \exists c_\alpha > 0 $ such that 
		\begin{equation}  \label{equn 14} 
			\bigwedge \{ t > 0 :  N_1 (  \beta_1 x_1  + \cdots + \beta_n x_n, K t   )  
			\geq 1 - \alpha  \} \geq \frac{c_\alpha}{\phi(s)}
		\end{equation}
		From (\ref{equn 13}) and (\ref{equn 14}), we have 
		\begin{align*}
			&	\frac{	n N_\alpha}{c_\alpha} \bigwedge \{ t > 0 :  N_1 ( x, K t )  \geq 1 - \alpha \} \geq
			\bigwedge \{ t > 0 :  N_2 ( Tx, t )  \geq \alpha * \cdots * \alpha \} \\
			\text{that is} ~ &  \bigwedge \{ t > 0 :  N_2 ( Tx, t )  \geq \alpha * \cdots * \alpha \} \leq M_\alpha' \bigwedge \{ t > 0 :  N_1 ( x, K t )  \geq 1 - \alpha \} ~ where ~ M_\alpha' = \frac{	n N_\alpha}{c_\alpha}  \\
			\text{that is} ~ & K \bigwedge \{ \frac{t}{K} > 0 :  N_2 ( Tx, K \cdot \frac{t}{K} )  \geq \alpha * \cdots * \alpha \} \leq \frac{M_\alpha'}{K} \bigwedge \{ K t > 0 :  N_1 ( x, K t )  \geq 1 - \alpha \} \\
			\text{that is} ~ & K \bigwedge \{ t > 0 :  N_2 ( Tx, K t )  \geq \alpha * \cdots * \alpha \} \leq \frac{M_\alpha'}{K} \bigwedge \{ t > 0 :  N_1 ( x, 
			t )  \geq 1 - \alpha \}. 
		\end{align*}
		Therefore finally we obtain, 
		\begin{equation} \label{equn 15}
			\bigwedge \{ t > 0 :  N_2 ( Tx, K t )  \geq \alpha * \cdots * \alpha \} \leq  M_\alpha \bigwedge \{ t > 0 :  N_1 ( x, 
			t )  \geq 1 - \alpha \} 
		\end{equation}
		where $ M_\alpha = \frac{M_\alpha '}{K^2} $. \\
		Since $ \alpha \geq \alpha * \cdots * \alpha $, so from (\ref{equn 15}) we have 
		\begin{align*}
		& \bigwedge \{ t > 0 :  N_2 ( Tx, K t )  \geq \alpha * \alpha * \cdots * \alpha \} \leq  M_\alpha \bigwedge \{ t > 0 :  N_1 ( x, 	t )  \geq 1 - \alpha * \alpha * \cdots * \alpha \}  ~~ \forall x \in X \\
		or ~ & 	\bigwedge \{ t > 0 :  N_2 ( Tx, K t )  \geq \beta( \alpha )  \} \leq  M_\alpha \bigwedge \{ t > 0 :  N_1 ( x, 
		t )  \geq 1 - \beta( \alpha ) \}  ~~ ~ \forall x \in X.
		\end{align*} 
		where $ \beta( \alpha ) = \alpha * \alpha * \cdots * \alpha $. \\
		Since $ \alpha \in (0, 1) $ is arbitrary, thus $T$ is fuzzy bounded.	
	\end{proof}	
\end{thm}	
\section{Operator fuzzy norm in $ BF(X, Y)$ }	
In this section, we define operator fuzzy norm in $ BF(X, Y)$ and finally study the completeness of $ BF(X, Y)$.
\begin{thm}  \label{thm 4}
	Let $ (X, N_1, \phi, K,   *_1) $ and $  (Y, N_2, \phi, K,   *_2) $ be two fuzzy strong $\phi$-b-normed linear spaces and $ *_2$ be lower semi-continuous. Let $ BF(X, Y ) $ denotes
	the set of all fuzzy bounded linear operators defined from $X$ to $Y$. 
	Then the mapping $ N : BF(X, Y) \times \mathbb{R} \rightarrow [0,1] $ defined by 
	$$ 	N ( T,s ) = \begin{cases}
	  \bigvee \{ \alpha \in (0,1) : \underset{x \in X\setminus \{ \theta\} }{\bigvee} \bigwedge \{ \frac{t}{d_{1 - \alpha}} > 0 : N_2(Tx, t )\geq \alpha \} \leq s  \} ~ ~ ~ for ~
	  (T,s ) \neq (0,0) \\
	  0 ~~ \hskip 238pt ~ for ~ ~ (T,s ) = (0,0)
	\end{cases} $$
is a fuzzy norm in $ BF(X, Y)$ with respect to underlying t-norm $ *_2 $.
\end{thm}	
\begin{proof}
	First we show that for  $ T \in  BF(X, Y)$,  $ \underset{x \in X\setminus \{ \theta\} }{\bigvee} \bigwedge \{ \frac{t}{d_{1 - \alpha}} > 0 : N_2(Tx, t )\geq \alpha \} $ 
	exists for each $  \alpha \in (0,1) $ and monotonically increasing  with respect to $  \alpha $. \\
	Since $ T \in  BF(X, Y)$, for each $  \alpha \in (0,1) $,  $ \exists M_\alpha > 0 $ such that 
	\begin{align*}
	&  \bigwedge \{ t > 0 : N_2 (  Tx, K t ) \geq \alpha \}  \leq  M_\alpha \bigwedge \{ t > 0 : N_1 (  x, t ) \geq 1 - \alpha \} ~ ~ 	\forall x \in X \\	
	\implies &  \bigwedge \{ t > 0 : N_2 (  Tx, K t ) \geq \alpha \}  \leq  M_\alpha d_{1 - \alpha} ~ ~ 	\forall x \in X \\
		\implies &  \bigwedge \{ \frac{t}{d_{1 - \alpha}} > 0 : N_2 (  Tx, K t ) \geq \alpha \}  \leq  M_\alpha ~ ~ 	\forall x(\neq \theta ) \in X. 
	\end{align*}
Let $  K t  = r $. Then we get, 
\begin{align*}
	&  \bigwedge \{ \frac{r}{Kd_{1 - \alpha}} > 0 : N_2 (  Tx, r ) \geq \alpha \}  \leq  M_\alpha ~ ~ 	\forall x(\neq \theta ) \in X \\
	\implies & \bigwedge \{ \frac{r}{d_{1 - \alpha}} > 0 : N_2 (  Tx, r ) \geq \alpha \}  \leq  K M_\alpha ~ ~ 	\forall x(\neq \theta ) \in X.
\end{align*}
Thus $  \underset{x \in X\setminus \{ \theta\} }{\bigvee} \bigwedge \{ \frac{t}{d_{1 - \alpha}} > 0 : N_2(Tx, t )\geq \alpha \} $ 
exists for each $  \alpha \in (0,1) $. \\
Next part, that is  monotonically increasing  with respect to $  \alpha $ is same as the proof in Theorem 4.1 \cite{6}. \\
Now we verify that $N$ satisfies (NI)-(NV). 
\begin{enumerate}[(i)]
	\item 
	Proof of (NI) and (NII) are same as in Theorem 4.1 \cite{6}.
	\item 
	For any scaler $ \lambda > 0 $, we have 
	\begin{align*}
		N( \lambda T, s )  = &   \bigvee \{ \alpha \in (0,1) : \underset{x \in X\setminus \{ \theta\} }{\bigvee} \bigwedge \{ \frac{t}{d_{1 - \alpha}} > 0 : 
		N_2((\lambda T) x, t )\geq 	\alpha \} \leq s  \}  \\
		= &   \bigvee \{ \alpha \in (0,1) : \underset{x \in X\setminus \{ \theta\} }{\bigvee} \bigwedge \{ \frac{t}{d_{1 - \alpha}} > 0 : 
		N_2(T x, \frac{t}{\phi(\lambda)} )\geq 	\alpha \} \leq s  \}  \\ 
		= &   \bigvee \{ \alpha \in (0,1) : \phi(\lambda) \underset{x \in X\setminus \{ \theta\} }{\bigvee} \bigwedge \{ \frac{t}{d_{1 - \alpha}} > 0 : 
		N_2(T x, t )\geq 	\alpha \} \leq s  \}  \\
		= &   \bigvee \{ \alpha \in (0,1) : \underset{x \in X\setminus \{ \theta\} }{\bigvee} \bigwedge \{ \frac{t}{d_{1 - \alpha}} > 0 : 
		N_2(T x, t )\geq 	\alpha \} \leq \frac{s}{\phi(\lambda) }  \}  \\
		= & N (T, \frac{s}{\phi(\lambda) } )
	\end{align*}
  So, (NIII) holds.
  \item
  We have to show that 
  $$ N (T_1 + T_2, s+ Kt ) \geq N (T_1, s) *_2 N (T_2, t) ~ ~ \forall s, t \in \mathbb{R}.  $$
  If possible suppose that the above relation does not hold. So $ \exists s_0, t_0 \in \mathbb{R}  $ such that \\
   $$ N (T_1 + T_2, s_0 + K t_0 ) < N (T_1, s_0) *_2 N (T_2, t_0). $$
   Choose $ \alpha_0 \in (0,1) $ such that 
   \begin{equation}  \label{equn 6}
   	 N (T_1 + T_2, s_0 + K t_0 ) < \alpha_0 < N (T_1, s_0) *_2 N (T_2, t_0). 
   \end{equation}
Since $ *_2 $ is lower semi-continuous,  $ \exists ~ \alpha_1, \alpha_2 \in (0,1)   $ where $  N (T_1, s_0) > \alpha_1 $ and   $ N (T_2, t_0) > \alpha_2 $ such that  
$ \alpha_1  *_2 \alpha_2 > \alpha_0 $. Now, 
\begin{align*}
	&  N (T_1, s_0) > \alpha_1  \\
	\implies & \underset{x \in X\setminus \{ \theta\} }{\bigvee} \bigwedge \{ \frac{s}{d_{1 - \alpha_1}} > 0 : 
	N_2(T_1 x, s )\geq 	\alpha_1 \} \leq s_0 \\
	\implies & \bigwedge \{ \frac{s}{d_{1 - \alpha_1}} > 0 : 	N_2(T_1 x, s )\geq 	\alpha_1 \} \leq s_0 ~ ~ 	\forall x(\neq \theta ) \in X. 
\end{align*}
Similarly,
\begin{align*}
	& K  \bigwedge \{ \frac{t}{d_{1 - \alpha_2}} > 0 : 	N_2(T_2 x, t )\geq 	\alpha_2 \} \leq K t_0 ~ ~ 	\forall x(\neq \theta ) \in X  \\
		\implies &  \bigwedge \{ \frac{K t}{d_{1 - \alpha_2}} > 0 : 	N_2(T_2 x, \frac{K t}{K} )\geq 	\alpha_2 \} \leq K t_0 ~ ~ 	\forall x(\neq \theta ) \in X  \\
		\implies &  \bigwedge \{ \frac{t}{d_{1 - \alpha_2}} > 0 : 	N_2(T_2 x, \frac{t}{K} )\geq 	\alpha_2 \} \leq K t_0 ~ ~ 	\forall x(\neq \theta ) \in X.  	
\end{align*}
Since $ \alpha_1 \geq \alpha_1 *_2 \alpha_2 > \alpha_0 $, so $ 1 - \alpha_0 > 1 - \alpha_1 $. Similarly, $ 1 - \alpha_0 > 1 - \alpha_2 $. Thus we get, 
$$ \bigwedge \{ \frac{s+t}{d_{1 - \alpha_0}} > 0 : 	N_2(T_1 x, s )\geq 	\alpha_1, ~ N_2(T_2 x, \frac{t}{K} )\geq 	\alpha_2 \} \leq s_0 + K t_0 
~ ~ 	\forall x(\neq \theta ) \in X. $$
Therefore,
\begin{align*}
&   \bigwedge \{ \frac{s+t}{d_{1 - \alpha_0}} > 0 : 	N_2( ( T_1 + T_2) x, s + K \frac{t}{K} ) \geq	\alpha_1 *_2 \alpha_2 \} \leq s_0 + K t_0 
~ ~ 	\forall x(\neq \theta ) \in X 	\\
\implies &  \bigwedge \{ \frac{s+t}{d_{1 - \alpha_0}} > 0 : 	N_2( ( T_1 + T_2) x, s + t ) \geq	\alpha_1 *_2 \alpha_2 \} \leq s_0 + K t_0 
~ ~ 	\forall x(\neq \theta ) \in X \\
\implies & \bigwedge \{ \frac{s+t}{d_{1 - \alpha_0}} > 0 : 	N_2( ( T_1 + T_2) x, s + t ) \geq \alpha_0 \} \leq s_0 + K t_0 
~ ~ 	\forall x(\neq \theta ) \in X  \\
\implies & \underset{x \in X\setminus \{ \theta\} }{\bigvee} \bigwedge \{ \frac{s+t}{d_{1 - \alpha_0}} > 0 : 	N_2( ( T_1 + T_2) x, s + t ) \geq	\alpha_0 \} 
\leq s_0 + K t_0 \\
\implies & N(T_1 + T_2, s_0 + K t_0 ) \geq \alpha_0.
\end{align*}
This contradicts the relation (\ref{equn 6}). \\
Hence (NIV): $ N (T_1 + T_2, s+ Kt ) \geq N (T_1, s) *_2 N (T_2, t) ~ ~ \forall s, t \in \mathbb{R} $ holds. 
\item 
The proof of the condition (NV) is same as the proof of the (NV) in Theorem 4.1\cite{6}.
\end{enumerate} 
\end{proof}	
\begin{ppn}\label{ppn 2}
Let $ (X, N, \phi, K, *) $  be a  fuzzy strong $\phi$-b-normed linear spaces and $ * $ be lower semi-continuous. Then limit of every l-fuzzy convergent sequence in $X$ 
is unique.
\begin{proof}
Let $ \beta \in (0,1)$. By the 	lower semi-continuity of $ * $, $ \exists \alpha \in (0,1) $ such that $ (1 - \alpha ) * (1 - \alpha ) > (1 - \beta ) $. \\
Let $ \{ x_n \} $ be an l-fuzzy convergent sequence in $X$ which converges to two different limits $ x $ and $y$. So, 
$$ \underset{n\rightarrow \infty}{lim}  \bigwedge \{ t > 0 : N (x_n - x, t) > 1 -\alpha  \} = 0, ~ \forall \alpha \in (0,1) $$
and 
$$ \underset{n\rightarrow \infty}{lim}  \bigwedge \{ t > 0 : N (x_n - y, t) > 1 -\alpha  \} = 0, ~ \forall \alpha \in (0,1). $$
Then for $ \epsilon > 0 $, $ \exists $ positive integers $ N_1 ( \alpha, \epsilon ) $ and $ N_2 ( \alpha, \epsilon ) $ such that 
$$   \bigwedge \{ t > 0 : N (x_n - x, t) > 1 -\alpha  \} < \frac{\epsilon}{2 K }, ~ \forall n \geq N_1 ( \alpha, \epsilon )  $$
and 
$$   \bigwedge \{ t > 0 : N (x_n - y, t) > 1 -\alpha  \} < \frac{\epsilon}{2 }, ~ \forall n \geq N_2 ( \alpha, \epsilon ). $$ 
Let $ N_0 = \max \{ N_1, N_2 \} $. Then, 
 $$ N ( x_n - x, \frac{\epsilon}{2 K } )  > 1 -\alpha, ~ \forall n \geq N_0  ( \alpha, \epsilon )  $$
 and
 $$ N ( x_n - y, \frac{\epsilon}{2 } )  > 1 -\alpha, ~ \forall n \geq N_0 ( \alpha, \epsilon ).  $$
 Now, 
 \begin{align*}
 	N ( x - y, \epsilon ) = & 	N ( x_n - y - x_n + x,  \frac{\epsilon}{2 } + K \cdot  \frac{\epsilon}{2 K } ) \\
 	\geq & N ( x_n - y, \frac{\epsilon}{2 } ) * N (  x_n - x,     \frac{\epsilon}{2 K } ) \\
 	> & (1 - \alpha ) * (1 - \alpha ) > (1 - \beta ).
 \end{align*}
Since, $ \epsilon > 0 $ and $ \beta \in (0,1)$ are arbitrary, so it implies that
\begin{align*}
	&  N ( x-y, t ) > \alpha, ~~ \forall t> 0, ~ ~    \alpha \in (0,1) \\
	\implies &  N ( x-y, t ) ~ = 1, ~~ \forall t> 0  \\
	\implies & x - y = \theta \\
		\implies & x = y.
\end{align*}
This completes the proof.
\end{proof}	
\end{ppn}	
\begin{thm} \label{thm 5}
Let $ (X, N_1, \phi, K, *_1 ) $ and $  (Y, N_2, \phi, K, *_2) $ be two fuzzy strong $\phi$-b-normed linear spaces and $ *_2$ be lower semi-continuous. If 
$    (Y, N_2, \phi, K, *_2) $ is l-fuzzy complete then  $ BF(X, Y ) $ 	is also l-fuzzy complete with respect to the underlying t-norm $ *_2 $. 
\begin{proof}
	Let $ \{ T_n \} $ be an l-fuzzy Cauchy sequence in $ BF(X, Y ) $. \\
	Then  $ \{ T_n x \} $ is an l-fuzzy Cauchy sequence in $ Y  $ for each $ x \in X $(Proof is similar to the proof of Theorem 5.1\cite{6}). \\
	Since $ Y  $ is l-fuzzy complete, so for each $ x \in X $, $ \exists y \in Y $ such that  $ \underset{n\rightarrow \infty}{lim} T_n x = y $. \\
	Thus we can define a function $T$ given by $ \underset{n\rightarrow \infty}{lim} T_n x = T x $. \\
	So,
	\begin{equation}  \label{equn 7}
		\underset{n\rightarrow \infty}{lim}  \bigwedge \{ t > 0 : N_2 (T_n x - T x, t) > 1 -\alpha  \} = 0  ~ ~~ \forall  x \in X, ~ ~ \forall \alpha \in (0,1). 
	\end{equation}
 First we show that $T$ is linear. \\
Now (\ref{equn 7}) implies that, for $ \epsilon > 0, ~ \exists N_1 ( \alpha, \epsilon ) $ such that 
\begin{align*}
	 &  ~ \bigwedge \{ t > 0 : N_2 (T_n x - T x, t) > 1 -\alpha  \} < \frac{\epsilon}{2} ~ ~ 
	 ~ \forall  x \in X,  \forall n \geq N_1 ( \alpha, \epsilon ) \\
	 \implies & N_2 (T_n x - T x, \frac{\epsilon}{2} ) > 1 -\alpha ~   ~~~ \forall  x \in X, ~ \forall n \geq N_1 ( \alpha, \epsilon ).  
\end{align*}
Similarly, 
$$ N_2 (T_n y - T y, \frac{\epsilon}{2 K } ) > 1 -\alpha ~   ~~~ \forall  y \in X, ~ \forall n \geq N_2 ( \alpha, \epsilon ). $$
Thus, 
$$ N_2 (T_n x - T x, \frac{\epsilon}{2} ) *_2  N_2 (T_n y - T y, \frac{\epsilon}{2 K } ) > ( 1 -\alpha) *_2  ( 1 -\alpha), ~~~~ \forall n \geq N_0 ( \alpha, \epsilon ) $$
where $ N_0 = \max \{ N_1, N_2 \} $. \\
So we get, 
$$ N_2 ( T_n (x+y) - ( Tx+Ty ), \frac{\epsilon}{2} + K \cdot \frac{\epsilon}{2K} ) \geq  ( 1 -\alpha) *_2  ( 1 -\alpha) ~~~~ \forall n \geq N_0 ( \alpha, \epsilon ). $$
Let $ \beta \in (0,1) $. Then $ \exists \alpha = \alpha(\beta) \in (0,1) $ such that $ ( 1 -\alpha) *_2  ( 1 -\alpha)  > (1-\beta) $. \\
Thus we have,
\begin{align*}
	& N_2 ( T_n (x+y) - ( Tx+Ty ), \epsilon ) > ( 1 -\beta)  ~~~~ \forall n \geq N_0  \\
	\implies &  \bigwedge \{ t > 0 : N_2 ( T_n (x+y) - ( Tx+Ty ), t) > 1 - \beta  \} \leq \epsilon ~ ~ \forall n \geq N_0  \\ 
	\implies &  \underset{n\rightarrow \infty}{lim} \bigwedge \{ t > 0 : N_2 ( T_n (x+y) - ( Tx+Ty ), t) > 1 - \beta  \} = 0.  
\end{align*} 
Since $ \beta \in (0,1)  $ is arbitrary, it follows that $ \{ T_n (x+y) \} $ is l-fuzzy convergent and converges to $ ( Tx+Ty ) $. So,
$$ \underset{n\rightarrow \infty}{lim} T_n (x+y) =  Tx + Ty  ~ ~ \text{ that is} ~ ~ T (x+y) =  Tx + Ty. $$	
Again for any scalar $ \lambda $ we have,
\begin{align*}
	& \underset{n\rightarrow \infty}{lim} T_n ( \lambda x ) = T ( \lambda x ) ~ ~~ \forall x \in X \\
	\implies & \underset{n\rightarrow \infty}{lim} \{ \lambda T_n  x \} = T ( \lambda x ) ~ ~~ \forall x \in X \\
	\implies & \lambda \underset{n\rightarrow \infty}{lim} T_n  x = T ( \lambda x ) ~ ~~ \forall x \in X \\
	\implies & T ( \lambda x )  = \lambda  ( T x ) ~ ~~ \forall x \in X. 
\end{align*}
Hence $T$ is linear. \\
Now we show that $T$ is fuzzy bounded. \\
Let $ \gamma \in (0,1) $.  By the lower semi-continuity of $ *_2 $, $ \exists \alpha = \alpha(\gamma) \in (0,1) $ such that  
$ ( 1 -\alpha) *_2  ( 1 -\alpha)  > \gamma $. \\
Now $ \{ T_n \} $ is an l-fuzzy Cauchy sequence in $ BF(X, Y ) $, so 
$$ \underset{m, n\rightarrow \infty}{lim}  \bigwedge \{ t > 0 : N_2 (T_m - T_n, t) > 1 -\alpha  \} = 0, ~ \forall \alpha \in (0,1). $$
Thus for a given $ \epsilon > 0 $ and for $ \alpha \in (0,1) $, there exists  positive integer $ N ( \alpha, \epsilon ) $  such that 
\begin{align*}
	& \bigwedge \{ t > 0 : N (T_m - T_n, t) > 1 -\alpha  \} < \frac{\epsilon}{2} ~ ~ ~ \forall m, n \geq N ( \alpha, \epsilon ) \\
	\implies &  N (T_m - T_n, \frac{\epsilon}{2} ) > 1 -\alpha ~ ~ ~ \forall m, n \geq N ( \alpha, \epsilon )
\end{align*}
which implies 
\begin{equation} \label{equn 8}
	\underset{x \in X\setminus \{ \theta\} }{\bigvee} \bigwedge \{ \frac{s}{d_{ \alpha}} > 0 : N_2 (T_n x - T_m x, s ) \geq 1 - \alpha \} \leq \frac{\epsilon}{2} ~
	~ ~ \forall m, n \geq N ( \alpha, \epsilon ).
\end{equation}	
From (\ref{equn 8}), for $ \alpha = \alpha(\gamma) \in (0,1), ~ \epsilon > 0, ~ \exists N'( \alpha(\gamma), \epsilon ) \in \mathbb{N} $ such that 
\begin{align*}
	& \bigwedge \{ \frac{s}{d_{ \alpha}} > 0 : N_2 (T_n x - T_m x, s ) \geq 1 - \alpha \} < \frac{\epsilon}{2} ~ ~ ~ \forall m, n \geq N' ( \alpha(\gamma), \epsilon ), ~~ x(\neq \theta ) \in X \\
	\implies & \bigwedge \{ s > 0 : N_2 (T_n x - T_m x, s d_{ \alpha} ) \geq 1 - \alpha \} < \frac{\epsilon}{2} ~ ~ ~ \forall m, n \geq N' ( \alpha(\gamma), \epsilon ), ~~ x(\neq \theta ) \in X 
\end{align*}
$$  $$
which implies 
\begin{equation} \label{equn 9}
 N_2 (T_n x - T_m x, \dfrac{d_\alpha \epsilon}{2} ) \geq 1 - \alpha, ~ ~ \forall m, n \geq N' ( \alpha(\gamma), \epsilon, x ) ~ ~~ x(\neq \theta ) \in X. 	
\end{equation}  
From (\ref{equn 7}), it follows that,   for $ \alpha = \alpha(\gamma) \in (0,1), ~ \epsilon > 0, ~ x(\neq \theta ) \in X, ~ \exists N''( \alpha(\gamma), \epsilon, x ) 
\in \mathbb{N} $ such that 
$$ \bigwedge \{ t > 0 : N_2 (T_m x - T x, t) > 1 -\alpha  \} < \dfrac{d_\alpha \epsilon}{2 K} ~ ~ ~ \forall m \geq N'' (\alpha(\gamma), \epsilon, x ).  $$
Therefore, 
\begin{equation} \label{equn 10}
	N_2 (T_m x - T x, \dfrac{d_\alpha \epsilon}{2 K} ) > 1 -\alpha ~ ~ ~ \forall m \geq N'' (\alpha(\gamma), \epsilon, x )
\end{equation}	   
Let $ N_0  (\alpha(\gamma), \epsilon, x ) = \max \{ N' ( \alpha(\gamma), \epsilon, x ), N'' (\alpha(\gamma), \epsilon, x )  \}    $. Then from (\ref{equn 9}) and (\ref{equn 10}) we 
have, 
\begin{align*}
	N_2 (T_n x - T x, d_\alpha \epsilon ) 	= & N_2 (T_n x - T_m x + T_m x - T x, \dfrac{d_\alpha \epsilon}{2 } + K \cdot \dfrac{d_\alpha \epsilon}{2 K} ) \\
	\geq & N_2 (T_n x - T_m x, \dfrac{d_\alpha \epsilon}{2} ) *_2 N_2 (T_m x - T x, \dfrac{d_\alpha \epsilon}{2 K} ) \\
	\geq & ( 1 -\alpha ) *_2 (1 -\alpha) \geq \gamma ~ ~ ~ \forall n \geq N_0 (\alpha(\gamma), \epsilon ), ~~ \forall  x(\neq \theta ) \in X.
\end{align*}
which implies  
\begin{equation} \label{equn 10.5}
	\bigwedge \{ t > 0 : 	N_2 (T_n x - T x,  t ) \geq \gamma \} \leq \epsilon d_\alpha  \leq \epsilon d_{1 - \gamma} ~ 
	~~ \forall n \geq N_0 (\alpha(\gamma), \epsilon ).
\end{equation}
Hence 
\begin{equation}  \label{equn 11}
 \bigwedge \{ t > 0 : 	N_2 (T_n x - T x, K t ) \geq \gamma \} \leq \frac{\epsilon}{K} 
\bigwedge \{ s > 0 : N_1 ( x, s ) \geq 1 - \gamma \} ~~ ~~ \forall n \geq N_0 (\alpha(\gamma), \epsilon ). 	
\end{equation}
This shows that $ T_n -T $ is fuzzy bounded $ \forall n \geq N_0 (\alpha(\gamma), \epsilon )$. \\
Since $ BF (X, Y ) $ is a linear space, so $ T = ( T -T_n ) + T_n \in BF (X, Y ) $. Thus $T$ is fuzzy bounded. \\
Again from (\ref{equn 10.5}), we get 
\begin{align*}
& \underset{x \in X\setminus \{ \theta\} }{\bigvee} \bigwedge \{ \frac{t}{d_{ 1 - \gamma}} > 0 : N_2 (T_n x - T x, t ) \geq 
\gamma \} \leq \epsilon ~ ~ ~ \forall m, n \geq N_0 (\alpha(\gamma), \epsilon ) \\
\implies  &  N (T_n - T, \epsilon ) \geq \gamma ~  ~ ~ \forall  n \geq N_0 (\alpha(\gamma), \epsilon ) \\
\implies  & \bigwedge \{ s > 0 : N (T_n - T, s ) \geq \gamma \} \leq \epsilon ~ ~ ~ \forall  n \geq N_0 (\alpha(\gamma), \epsilon ) 
\end{align*} 
Since $ \epsilon > 0 $ is arbitrary, so from above we have, 
$$  \underset{n\rightarrow \infty}{lim} \bigwedge \{ s > 0 : N (T_n - T, s ) \geq \gamma \} = 0. $$ 

Since $ \gamma \in (0,1) $ is arbitrary, so it follows that 
$$ \underset{n\rightarrow \infty}{lim} \bigwedge \{ s > 0 : N (T_n - T, s ) \geq \alpha \} = 0 ~~~ \forall \alpha \in (0,1).  $$
So $ \{ T_n \} $ is l-fuzzy convergent and converges to $T \in BF (X, Y ) $. Thus  $ BF (X, Y ) $ is l-fuzzy complete.
\end{proof}
\end{thm}
%
%
%
%
%
~\\
%
%
\textbf{Conclusion:}  
In this paper, we are able to extend some concepts of fuzzy normed linear spaces to fuzzy strong $\phi$-b-normed linear spaces. Definitions of fuzzy bounded linear operator and fuzzy continuous operator are given and their relation is studied. Idea of operator norm in BF(X,Y) is introduced and completeness of BF(X,Y) is established.  \\
 There are huge  scope for development in continuation with the  results  of this manuscript such as four fundamental theorems of functional analysis, uniform convex and strictly  convex spaces in  fuzzy setting. 
 \\
~\\
\textbf{Acknowledgment:}  The author AD is grateful to University Grant Commission (UGC), 
New Delhi, India for awarding her senior research fellowship [Grant No.1221/(CSIRNETJUNE2019)].
The authors are  thankful  to the Department of Mathematics, Siksha-Bhavana, Visva-Bharati.   \\

\end{document}